\documentclass[12pt]{article}%
\usepackage{amsfonts}
\usepackage{sw20bams}
\usepackage{amsmath}
\usepackage{amssymb}
\usepackage{graphicx}%
\providecommand{\U}[1]{\protect\rule{.1in}{.1in}}
\begin{document}

\title{Exercises in Iterational Asymptotics}
\author{Steven Finch}
\date{March 6, 2025}
\maketitle

\begin{abstract}
The problems and solutions contained here, all associated with nonlinear
recurrences and long-term trends, are new (as far as is known). \ 

\end{abstract}

\footnotetext{Copyright \copyright \ 2025 by Steven R. Finch. All rights
reserved.}

\section{Premi\`{e}re exercice}

\textbf{Consider initially} the famous recurrence%
\[%
\begin{array}
[c]{ccccc}%
x_{k}=p\,x_{k-1}\left(  1-x_{k-1}\right)  &  & \text{for }k\geq1\text{;} &  &
0<x_{0}<1
\end{array}
\]
where $0<p<1$. \ Quantify the convergence rate of $x_{k}$ as $k\rightarrow
\infty$.

Clearly $0<x_{k}<1$ and%
\[
x_{k}<p\,x_{k-1}<p^{2}x_{k-2}<p^{3}x_{k-3}<p^{4}x_{k-4}%
\]
thus $x_{k}<p^{k}x_{0}$ for all $k$. \ Observe that%
\begin{align*}
x_{k}  &  =p\,x_{k-1}\left(  1-x_{k-1}\right) \\
&  =p^{2}x_{k-2}\left(  1-x_{k-2}\right)  \left(  1-x_{k-1}\right) \\
&  =p^{3}x_{k-3}\left(  1-x_{k-3}\right)  \left(  1-x_{k-2}\right)  \left(
1-x_{k-1}\right) \\
&  =p^{k}x_{0}%
{\displaystyle\prod\limits_{j=0}^{k-1}}
\left(  1-x_{j}\right)
\end{align*}
hence%
\[
C=\lim_{k\rightarrow\infty}\frac{x_{k}}{p^{k}}=x_{0}%
{\displaystyle\prod\limits_{j=0}^{\infty}}
\left(  1-x_{j}\right)
\]
exists and is nonzero since%
\[%
{\displaystyle\sum\limits_{j=0}^{\infty}}
x_{j}<x_{0}%
{\displaystyle\sum\limits_{j=0}^{\infty}}
p^{j}%
\]
converges.

Table 1. \textit{Numerical estimates of }$C$\textit{, given }$x_{0}=1/2$
\textit{and selected values of }$p>0$: \textit{no closed-form expressions are
known}\
\[%
\begin{tabular}
[c]{|l|l|l|l|l|}\hline
$p$ & $C$ &  & $p$ & $C$\\\hline
1/5 & $0.234690787230465...$ &  & 3/5 & $0.176983588618567...$\\\hline
1/4 & $0.229832778573153...$ &  & 2/3 & $0.161059687971223...$\\\hline
1/3 & $0.220577540168322...$ &  & 3/4 & $0.136649472578135...$\\\hline
2/5 & $0.211947268934865...$ &  & 4/5 & $0.118823329484862...$\\\hline
1/2 & $0.196453426377889...$ &  & 1 & $1.767993786136154...^{\ast}$\\\hline
\end{tabular}
\ \ \ \ \ \ \
\]
The entry corresponding to $p=1$ is starred \cite{S1-exercis} because $C$ is
defined differently than for $p<1$:%
\[
C=-\lim_{k\rightarrow\infty}k^{2}\left(  x_{k}-\frac{1}{k}+\frac{\ln(k)}%
{k^{2}}\right)
\]
and details are found in \cite{F1-exercis, F4-exercis, Sc-exercis}.

\textbf{Consider finally }the less-famous recurrence%
\[%
\begin{array}
[c]{ccccc}%
x_{k}=p\,x_{k-1}\left(  1+x_{k-1}\right)  &  & \text{for }k\geq1\text{;} &  &
0<x_{0}<\dfrac{1-p}{p}%
\end{array}
\]
where $0<p<1$. \ Again, quantify the convergence rate of $x_{k}$ as
$k\rightarrow\infty$.

Replacing $1-x_{k-1}$ by $1+x_{k-1}$ renders the task of bounding $x_{k}$ more
subtle. \ Note that $(1-p)/p$ is the only nonzero fixed point of
$f(x)=p\,x\left(  1+x\right)  $:%
\[%
\begin{array}
[c]{ccccc}%
1=p(1+x) &  & \text{when} &  & 1-p=p\,x.
\end{array}
\]
If $x_{0}=(1-p-\varepsilon)/p$ for some $0<\varepsilon<1-p$, then%
\[
p(1+x_{0})=p+(1-p-\varepsilon)=1-\varepsilon<1
\]
hence%
\[
x_{1}=p\,x_{0}\left(  1+x_{0}\right)  =p\left(  1+x_{0}\right)  x_{0}%
<(1-\varepsilon)x_{0}<x_{0}.
\]
A more geometric approach involves graphing the parabola $y=f(x)$ and the
diagonal $y=x$:\ since $f^{\prime}(0)=p<1$, the curve dips below the line at
the left endpoint and, by continuity, does not cross it again until the right
endpoint. \ More generally, $x_{k}<x_{k-1}$ for all $k$. \ It follows that%
\begin{align*}
x_{k}  &  <p\,x_{k-1}\left(  1+x_{0}\right)  =p\left(  1+x_{0}\right)
x_{k-1}=\left(  1-\varepsilon\right)  x_{k-1}\\
&  <\left(  1-\varepsilon\right)  ^{2}x_{k-2}<\left(  1-\varepsilon\right)
^{3}x_{k-3}<\left(  1-\varepsilon\right)  ^{4}x_{k-4}%
\end{align*}
thus $x_{k}<\left(  1-\varepsilon\right)  ^{k}x_{0}$ for all $k$. \ Observe
that%
\begin{align*}
x_{k}  &  =p\,x_{k-1}\left(  1+x_{k-1}\right) \\
&  =p^{2}x_{k-2}\left(  1+x_{k-2}\right)  \left(  1+x_{k-1}\right) \\
&  =p^{3}x_{k-3}\left(  1+x_{k-3}\right)  \left(  1+x_{k-2}\right)  \left(
1+x_{k-1}\right) \\
&  =p^{k}x_{0}%
{\displaystyle\prod\limits_{j=0}^{k-1}}
\left(  1+x_{j}\right)
\end{align*}
hence%
\[
C=\lim_{k\rightarrow\infty}\frac{x_{k}}{p^{k}}=x_{0}%
{\displaystyle\prod\limits_{j=0}^{\infty}}
\left(  1+x_{j}\right)
\]
exists and is nonzero since%
\[%
{\displaystyle\sum\limits_{j=0}^{\infty}}
x_{j}<x_{0}%
{\displaystyle\sum\limits_{j=0}^{\infty}}
(1-\varepsilon)^{j}%
\]
converges.

Table 2. \textit{Numerical estimates of }$C$\textit{, given }$x_{0}%
=(1-p)/(2p)$: \textit{no closed-form expressions are known}\
\[%
\begin{tabular}
[c]{|l|l|l|l|l|}\hline
$p$ & $C$ &  & $p$ & $C$\\\hline
1/5 & $24.539007835941751...$ &  & 3/5 & $1.015970842139591...$\\\hline
1/4 & $13.119009853937092...$ &  & 2/3 & $0.690744393761287...$\\\hline
1/3 & $5.896477923507413...$ &  & 3/4 & $0.415551960439528...$\\\hline
2/5 & $3.529895194705441...$ &  & 4/5 & $0.295525160728184...$\\\hline
1/2 & $1.832010583354543...$ &  & 1 & $1.597910218031873...^{\ast}$\\\hline
\end{tabular}
\ \ \ \ \ \ \ \ \
\]
The entry corresponding to $p=1$ is starred \cite{S2-exercis} because again
$C$ is defined differently than for $p<1$:%
\[%
\begin{array}
[c]{ccccccc}%
C=\lim\limits_{k\rightarrow\infty}x_{k}^{(2^{-k})} &  & \text{where} &  &
x_{0}=1 &  & \text{(not }0\text{).}%
\end{array}
\]
This constant appears elsewhere, although thinly veiled \cite{AS-exercis,
F5-exercis}: \
\[
\sqrt{C}=\lim\limits_{k\rightarrow\infty}y_{k}^{(2^{-k-1})}%
=1.264084735305301...
\]
where $y_{k}=1+x_{k}$; we have
\[%
\begin{array}
[c]{ccccc}%
y_{k}=y_{k-1}^{2}-y_{k-1}+1 &  & \text{for }k\geq1\text{;} &  & y_{0}=2
\end{array}
\]
and the latter is known as Sylvester's sequence \cite{S3-exercis}. \ Also,
letting $z_{k}=1+2x_{k}$, we have
\[%
\begin{array}
[c]{ccccc}%
z_{k}=\dfrac{1}{2}\left(  z_{k-1}^{2}+1\right)   &  & \text{for }%
k\geq1\text{;} &  & z_{0}=3
\end{array}
\]
and the latter possesses an interesting connection with Pythagorean triples
\cite{S4-exercis}. \ Not only is $C$ irrational \cite{WZ-exercis}, it is also
transcendental \cite{Du-exercis}.

\section{Deuxi\`{e}me exercice}

\textbf{For two distinct} starting values $0<x_{0}<1$, determine numerically%
\[
C=C(x_{0})=-\frac{1}{8}\lim_{k\rightarrow\infty}k^{3}\left(  x_{k}-\frac
{4}{k^{2}}+12\frac{\ln(k)}{k^{3}}\right)
\]
where
\[%
\begin{array}
[c]{ccc}%
x_{k}=x_{k-1}\left(  1-\sqrt{x_{k-1}}\right)  &  & \text{for }k\geq1\text{.}%
\end{array}
\]
Using $C$, find the asymptotic expansion of $x_{k}$ to order $1/k^{6}$.

We shall accomplish the steps in reverse order. \ This work builds on
\cite{dB-exercis, BR-exercis, ML-exercis}. \ To conserve space, formulaic
knowledge of sections 1, 2, 3 of \cite{F4-exercis} is assumed. \ For the
function $f(x)=x\left(  1-\sqrt{x}\right)  $, we have $\tau=1/2$,%
\[
\{a_{m}\}_{m=1}^{7}=\left\{  -1,0,0,0,0,0,0\right\}
\]
and $\lambda=2$; consequently
\[
\{b_{j}\}_{j=1}^{6}=\left\{  \frac{3}{2},\frac{5}{2},\frac{35}{8},\frac{63}%
{8},\frac{231}{16},\frac{429}{16}\right\}  ,
\]%
\[
\{a_{0j}\}_{j=1}^{6}=\left\{  1,1,\frac{4}{3},2,\frac{16}{5},\frac{16}%
{3}\right\}  ,
\]%
\[
\{c_{i}\}_{i=1}^{5}=\left\{  1,\frac{15}{16},\frac{35}{24},\frac{14}{5}%
,\frac{448}{75}\right\}  ,
\]%
\[
T_{2}=\frac{3}{2}X-1,
\]%
\[
T_{3}=-\dfrac{3}{4}X^{2}+\dfrac{13}{4}X-\dfrac{39}{16},
\]%
\[
T_{4}=\dfrac{1}{2}X^{3}-\frac{35}{8}X^{2}+\dfrac{39}{4}X-\dfrac{587}{96},
\]

\[
T_{5}=-\frac{3}{8}X^{4}+\frac{41}{8}X^{3}-\frac{339}{16}X^{2}+\frac{1055}%
{32}X-\frac{5451}{320},
\]%
\[
T_{6}=\frac{3}{10}X^{5}-\frac{91}{16}X^{4}+\frac{575}{16}X^{3}-\frac{3127}%
{32}X^{2}+\frac{37629}{320}X-\frac{245957}{4800}%
\]
and
\[
P_{2}=\frac{3}{4}X^{2}+\frac{3}{2}X+2,
\]%
\[
P_{3}=\frac{1}{2}X^{3}+\dfrac{21}{8}X^{2}+\dfrac{25}{4}X+\dfrac{39}{8},
\]%
\[
P_{4}=\frac{5}{16}X^{4}+\dfrac{47}{16}X^{3}+\dfrac{47}{4}X^{2}+\dfrac{345}%
{16}X+\dfrac{731}{48},
\]%
\[
P_{5}=\frac{3}{16}X^{5}+\frac{171}{64}X^{4}+\dfrac{517}{32}X^{3}+\dfrac
{1599}{32}X^{2}+\dfrac{2497}{32}X+\dfrac{7791}{160},
\]%
\[
P_{6}=\frac{7}{64}X^{6}+\frac{1377}{640}X^{5}+\frac{4645}{256}X^{4}%
+\dfrac{2641}{32}X^{3}+\dfrac{27073}{128}X^{2}+\dfrac{11499}{40}%
X+\dfrac{3091081}{19200}.
\]
The remarkable formula connecting $P_{m}$ and asymptotics of $x_{k}%
=f(x_{k-1})$ is%
\[
x_{k}\sim\left(  \frac{\lambda}{k}\right)  ^{1/\tau}\left\{  1+%
{\displaystyle\sum\limits_{m=1}^{6}}
P_{m}\left(  -\frac{1}{\tau}\left[  b_{1}\ln(k)+C\right]  \right)  \frac
{1}{k^{m}}\right\}
\]
which implies%
\begin{align*}
x_{k}  &  \sim\frac{4}{k^{2}}-12\frac{\ln(k)}{k^{3}}-\frac{8C}{k^{3}}%
+27\frac{\ln(k)^{2}}{k^{4}}+\left(  36C-18\right)  \frac{\ln(k)}{k^{4}%
}+\left(  12C^{2}-12C+8\right)  \frac{1}{k^{4}}\\
&  -54\frac{\ln(k)^{3}}{k^{5}}-\left(  108C-\frac{189}{2}\right)  \frac
{\ln(k)^{2}}{k^{5}}-\left(  72C^{2}-126C+75\right)  \frac{\ln(k)}{k^{5}}\\
&  -\left(  16C^{3}-42C^{2}+50C-\frac{39}{2}\right)  \frac{1}{k^{5}}%
+\frac{405}{4}\frac{\ln(k)^{4}}{k^{6}}+\left(  270C-\frac{1269}{4}\right)
\frac{\ln(k)^{3}}{k^{6}}\\
&  +\left(  270C^{2}-\frac{1269}{2}C+423\right)  \frac{\ln(k)^{2}}{k^{6}%
}+\left(  120C^{3}-423C^{2}+564C-\frac{1035}{4}\right)  \frac{\ln(k)}{k^{6}}\\
&  +\left(  20C^{4}-94C^{3}+188C^{2}-\frac{345}{2}C+\frac{731}{12}\right)
\frac{1}{k^{6}}.
\end{align*}
Assuming $x_{0}=1/2$ (the midpoint), the constant $C$ is estimated to be%
\[
C=1.98803983644549695008812308629512...
\]
by a simple numerical method \cite{F1-exercis} using the preceding expansion.
\ Assuming $x_{0}=4/9$ (the argument at which $C$ is minimal), we have instead%
\[
C=1.96846882098495471088450855794395....
\]

A similar procedure applies to the recurrence%
\[%
\begin{array}
[c]{ccc}%
x_{k}=x_{k-1}\left(  1-x_{k-1}^{2}\right)  &  & \text{for }k\geq1.
\end{array}
\]
We did this earlier \cite{F2-exercis} but using a different approach; see also
\cite{LSZ-exercis}.\ \ Values of $C$ were found for $x_{0}=1/2$ and for
$x_{0}=1/\sqrt{3}$ (again, the argument that minimizes $C$). \ We leave
analyses of%
\[%
\begin{array}
[c]{ccccc}%
u_{k}=u_{k-1}\left(  1-\dfrac{1}{2}u_{k-1}^{2}\right)  , &  & v_{k}%
=v_{k-1}\cos(v_{k-1}), &  & w_{k}=w_{k-1}\exp\left(  -\dfrac{w_{k-1}^{2}}%
{2}\right)
\end{array}
\]
for an interested reader.

\section{Troisi\`{e}me exercice}

\textbf{For two distinct} parameter values $q>1$, determine numerically%
\[
C=C(q)=q\,\lim_{k\rightarrow\infty}k^{(q-1)/q}\left(  q^{(q-1)/q}%
\,x_{k}-q\,k^{1/q}-\frac{q-1}{2q}\frac{\ln(k)}{k^{(q-1)/q}}\right)
\]
where%
\[%
\begin{array}
[c]{ccccc}%
x_{k}=x_{k-1}+\dfrac{1}{x_{k-1}^{q-1}} &  & \text{for }k\geq1\text{;} &  &
x_{0}=1\text{.}%
\end{array}
\]
Using $C$, find the asymptotic expansion of $x_{k}$ to order $1/k^{(3q-1)/q}$.

Let $y_{k}=x_{k}^{q}$. \ From%
\[
\left(  \frac{y_{k}}{y_{k-1}}\right)  ^{1/q}=\frac{x_{k}}{x_{k-1}}=1+\dfrac
{1}{x_{k-1}^{q}}=1+\frac{1}{y_{k-1}}%
\]
we have%
\[
y_{k}=y_{k-1}\left(  1+\frac{1}{y_{k-1}}\right)  ^{q}=y_{k-1}\,\varphi\left(
\frac{1}{y_{k-1}}\right)
\]
where $\varphi(z)=(1+z)^{q}$. \ Note that $\varphi(z)>1$ for all $z>0$,
$\varphi(0)=1$ and the derivative $\varphi^{\prime}(z)$ satisfies
$\varphi^{\prime}(0)\neq0$; \ also%
\[
\alpha=\varphi^{\prime}(0)=q,
\]

\[
\beta=\frac{\varphi^{\prime\prime}(0)}{2\varphi^{\prime}(0)}=-\frac{1}%
{2}+\frac{q}{2},
\]%
\[
\gamma=\frac{\varphi^{\prime\prime}(0)^{2}}{4\varphi^{\prime}(0)^{3}}=\frac
{1}{4q}-\frac{1}{2}+\frac{q}{4},
\]%
\[
\delta=-\frac{\varphi^{\prime\prime}(0)}{4\varphi^{\prime}(0)}-\frac
{\varphi^{\prime\prime\prime}(0)-3C\,\varphi^{\prime\prime}(0)}{6\varphi
^{\prime}(0)^{2}}+\frac{\varphi^{\prime\prime}(0)^{2}}{4\varphi^{\prime
}(0)^{3}}=-\frac{1}{12q}-\frac{C}{2q}+\frac{1}{4}+\frac{C}{2}-\frac{q}{6}.
\]
By Theorem 6 of \cite{P1-exercis} and Theorem 5 of \cite{P2-exercis},%
\[
y_{k}\sim\alpha\,k+\beta\ln(k)+C+\gamma\frac{\ln(k)}{k}+\delta\frac{1}{k}%
\]
as $k\rightarrow\infty$ (beware: the lead coefficient $1/2$ of $\delta$ in
\cite{P2-exercis} should be $1/4$). \ Let $r=1/q$. \ By Proposition 7 of
\cite{P2-exercis},%
\begin{align*}
x_{k}  &  \sim\alpha^{r}k^{r}+\frac{r\,\beta}{\alpha^{1-r}}\frac{\ln
(k)}{k^{1-r}}+\frac{r\,C}{\alpha^{1-r}}\frac{1}{k^{1-r}}+\frac{r(r-1)\beta
^{2}}{2\alpha^{2-r}}\frac{\ln(k)^{2}}{k^{2-r}}\\
&  +\frac{r(r-1)\beta\,C+r\,\alpha\,\gamma}{\alpha^{2-r}}\frac{\ln(k)}%
{k^{2-r}}+\frac{r(r-1)C^{2}+2r\,\alpha\,\delta}{2\alpha^{2-r}}\frac{1}%
{k^{2-r}}%
\end{align*}
which, after multiplying both sides by $\alpha^{1-r}$, becomes%
\begin{align*}
q^{1-r}\,x_{k}  &  \sim q\,k^{r}+\left(  -\frac{1}{2q}+\frac{1}{2}\right)
\frac{\ln(k)}{k^{1-r}}+\frac{C}{q}\frac{1}{k^{^{1-r}}}-\left(  -\frac
{1}{8q^{3}}+\frac{3}{8q^{2}}-\frac{3}{8q}+\frac{1}{8}\right)  \frac{\ln
(k)^{2}}{k^{2-r}}\\
&  +\left[  \left(  \frac{1}{4q^{2}}-\frac{1}{2q}+\frac{1}{4}\right)  -\left(
\frac{1}{2q^{3}}-\frac{1}{q^{2}}+\frac{1}{2q}\right)  C\right]  \frac{\ln
(k)}{k^{2-r}}\\
&  +\left[  -\left(  \frac{1}{12q^{2}}-\frac{1}{4q}+\frac{1}{6}\right)
+\left(  -\frac{1}{2q^{2}}+\frac{1}{2q}\right)  C-\left(  -\frac{1}{2q^{3}%
}+\frac{1}{2q^{2}}\right)  C^{2}\right]  \frac{1}{k^{2-r}}.
\end{align*}
In particular,
\begin{align*}
2^{1/2}\,x_{k}  &  \sim2k^{1/2}+\frac{1}{4}\frac{\ln(k)}{k^{1/2}}+\frac{C}%
{2}\frac{1}{k^{1/2}}-\frac{1}{64}\frac{\ln(k)^{2}}{k^{3/2}}\\
&  +\left(  \frac{1}{16}-\frac{C}{16}\right)  \frac{\ln(k)}{k^{3/2}}+\left(
-\frac{1}{16}+\frac{C}{8}-\frac{C^{2}}{16}\right)  \frac{1}{k^{3/2}}%
\end{align*}
when $q=2$;%

\begin{align*}
3^{2/3}\,x_{k}  &  \sim3k^{1/3}+\frac{1}{3}\frac{\ln(k)}{k^{2/3}}+\frac{C}%
{3}\frac{1}{k^{2/3}}-\frac{1}{27}\frac{\ln(k)^{2}}{k^{5/3}}\\
&  +\left(  \frac{1}{9}-\frac{2C}{27}\right)  \frac{\ln(k)}{k^{5/3}}+\left(
-\frac{5}{54}+\frac{C}{9}-\frac{C^{2}}{27}\right)  \frac{1}{k^{5/3}}%
\end{align*}
when $q=3$;%

\begin{align*}
\left(  \frac{3}{2}\right)  ^{1/3}\,x_{k}  &  \sim\frac{3}{2}k^{2/3}+\frac
{1}{6}\frac{\ln(k)}{k^{1/3}}+\frac{2C}{3}\frac{1}{k^{1/3}}-\frac{1}{216}%
\frac{\ln(k)^{2}}{k^{4/3}}\\
&  +\left(  \frac{1}{36}-\frac{C}{27}\right)  \frac{\ln(k)}{k^{4/3}}+\left(
-\frac{1}{27}+\frac{C}{9}-\frac{2C^{2}}{27}\right)  \frac{1}{k^{4/3}}.
\end{align*}
when $q=3/2$. \ 

To obtain $C$ to high numerical precision requires more terms in the above
expansions. \ The elegant technique from \cite{BR-exercis, ML-exercis}, useful
in Section 2, does not apply here (since $x_{k}\nrightarrow0$). \ We turn to a
brute-force matching-coefficient method, demonstrated in \cite{F1-exercis,
Sc-exercis, F2-exercis, F3-exercis}. \ For $q=2$, we expand%
\[%
\begin{array}
[c]{ccccccc}%
(k+1)^{1/2}, &  & \dfrac{\ln(k+1)^{\ell}}{(k+1)^{1/2}}, &  & \dfrac
{\ln(k+1)^{m}}{(k+1)^{3/2}}, &  & \dfrac{\ln(k+1)^{n}}{(k+1)^{5/2}}%
\end{array}
\]
for $\ell=1,0$; $m=2,1,0$; $n=3,2,1,0$ and compare a series for $x_{k+1}$ with
a series for $x_{k}+x_{k}^{-1}$, yielding additional terms%
\begin{align*}
&  \frac{1}{512}\frac{\ln(k)^{3}}{k^{5/2}}+\left(  -\frac{1}{64}+\frac
{3C}{256}\right)  \frac{\ln(k)^{2}}{k^{5/2}}+\left(  \frac{5}{128}-\frac
{C}{16}+\frac{3C^{2}}{128}\right)  \frac{\ln(k)}{k^{5/2}}\\
&  +\left(  -\frac{11}{384}+\frac{5C}{64}-\frac{C^{2}}{16}+\frac{C^{3}}%
{64}\right)  \frac{1}{k^{5/2}}.
\end{align*}
For $q=3$, we expand%
\[%
\begin{array}
[c]{ccccccc}%
(k+1)^{1/3}, &  & \dfrac{\ln(k+1)^{\ell}}{(k+1)^{2/3}}, &  & \dfrac
{\ln(k+1)^{m}}{(k+1)^{5/3}}, &  & \dfrac{\ln(k+1)^{n}}{(k+1)^{8/3}}%
\end{array}
\]
for $\ell$, $m$, $n$ identical with before and compare a series for $x_{k+1}$
with a series for $x_{k}+x_{k}^{-2}$, yielding additional terms%
\begin{align*}
&  \frac{5}{729}\frac{\ln(k)^{3}}{k^{8/3}}+\left(  -\frac{7}{162}+\frac
{5C}{243}\right)  \frac{\ln(k)^{2}}{k^{8/3}}+\left(  \frac{43}{486}-\frac
{7C}{81}+\frac{5C^{2}}{243}\right)  \frac{\ln(k)}{k^{8/3}}\\
&  +\left(  -\frac{1}{18}+\frac{43C}{486}-\frac{7C^{2}}{162}+\frac{5C^{3}%
}{729}\right)  \frac{1}{k^{8/3}}.
\end{align*}
From such enhancements, letting $c=C/q$, we obtain%
\[
c(2)=0.8615711875687117305317813...=(0.6092228292047829402293060...)\sqrt{2},
\]%
\[
c(3)=1.3784186157718345713984647....
\]
Finding $c(3/2)=0.8010888849...$ accurately is more difficult, owing to the
square root in the iteration. \ Introducing a power series approximation for
it in the brute-force method might be an option; augmenting the theory in
\cite{P1-exercis, P2-exercis} is another. \ Most preferable, however, would be
a generalization of the algorithm presented in \cite{BR-exercis, ML-exercis}
to encompass certain cases where $x_{k}\rightarrow\infty$. \ The latter would
be a powerful broadening of available tools for asymptotic expansion.

\section{Acknowledgements}

Robert Israel and Anthony Quas gave the simple proof that $C$ exists for the
recurrence $x_{k}=(1-p)+p\,x_{k-1}^{2}$ when $0<p<1/2$; their technique was
extended in \cite{F1-exercis} to $1/2<p<1$ and applied to $x_{k}%
=\sqrt{1+x_{k-1}}$ in \cite{F3-exercis}. \ The creators of Mathematica earn my
gratitude every day:\ this paper could not have otherwise been written.

\section{Addendum}

A certain reciprocity has been found, permitting a calculation (left
unfinished in Section 3) to be completed:%
\[
c(3/2)=\Lambda=0.8010888849039666437110775...
\]
where%
\[
-\lim_{k\rightarrow\infty}k^{5/3}\left(  \left(  \frac{3}{2}\right)  ^{2/3}%
\xi_{k}-\dfrac{1}{k^{2/3}}+\dfrac{1}{9}\dfrac{\ln(k)}{k^{5/3}}\right)
=\dfrac{2\Lambda}{3},
\]%
\[%
\begin{array}
[c]{ccccc}%
\xi_{k}=\dfrac{\xi_{k-1}}{1+\xi_{k-1}^{3/2}} &  & \text{for }k\geq1\text{;} &
& \xi_{0}=1
\end{array}
\]
and this iteration (fortunately!) may be treated via \cite{BR-exercis,
ML-exercis}. \ In the same way,%
\[
-\lim_{k\rightarrow\infty}k^{4/3}\left(  3^{1/3}\eta_{k}-\dfrac{1}{k^{1/3}%
}+\dfrac{1}{9}\dfrac{\ln(k)}{k^{4/3}}\right)  =\dfrac{c(3)}{3},
\]%
\[
-\lim_{k\rightarrow\infty}k^{3/2}\left(  2^{1/2}\zeta_{k}-\dfrac{1}{k^{1/2}%
}+\dfrac{1}{8}\dfrac{\ln(k)}{k^{3/2}}\right)  =\dfrac{c(2)}{2}%
\]
where%
\[%
\begin{array}
[c]{cccccc}%
\eta_{k}=\dfrac{\eta_{k-1}}{1+\eta_{k-1}^{3}} & \text{and} & \zeta_{k}%
=\dfrac{\zeta_{k-1}}{1+\zeta_{k-1}^{2}} & \text{for }k\geq1\text{;} &  &
\eta_{0}=\zeta_{0}=1.
\end{array}
\]
We report on this phenomenon in greater depth \cite{F6-exercis} and provide
more exercises \cite{F7-exercis}.

\end{document}